\documentclass[runningheads]{svjour2}
\smartqed  

\usepackage{graphicx,subfigure}

\usepackage{amsfonts,amsmath,bm,hhline}

\DeclareGraphicsExtensions{.jpg,.png,.eps,.eps,.gif}

\journalname{Applicable Algebra in Engineering, Communication and Computing}
\begin{document}

\title{Algebraic change-point detection 
}

\titlerunning{Change-point detection}        

\author{Michel~Fliess~\and~C\'{e}dric~Join~\and~Mamadou~Mboup}

\authorrunning{Fliess, Join \& Mboup} 

\institute{M. Fliess \at INRIA-ALIEN \& LIX (CNRS, UMR 7161),
\'Ecole polytechnique,
91128 Palaiseau, France \\
              \email{Michel.Fliess@polytechnique.edu}           
           \and
           C. Join \at
           INRIA-ALIEN \& CRAN (CNRS, UMR 7039),
           Nancy-Universit\'{e}, BP 239, \\
           54506 Vand\oe{}uvre-l\`{e}s-Nancy, France \\
           \email{Cedric.Join@cran.uhp-nancy.fr}
             \and
             M. Mboup \at
             INRIA-ALIEN \& CReSTIC, UFR Sciences Exactes et
             Naturelles, \\
Universit\'{e} de Reims Champagne-Ardenne, Moulin de la Housse, B.P.
1039, \\ 51687 Reims cedex 2, France
\\  \email{Mamadou.Mboup@univ-reims.fr} }

\date{Received: date / Revised: date}

\maketitle

\begin{abstract}
Elementary techniques from operational calculus, differential
algebra, and noncommutative algebra lead to a new approach for
change-point detection, which is an important field of investigation
in various areas of applied sciences and engineering. Several
successful numerical experiments are presented.
\keywords{Change-point detection \and Identifiability \and
Operational calculus \and Differential algebra \and Noncommutative
algebra \and Holonomic functions}
\end{abstract}

\section{Introduction}\label{intro}
Let $f: \mathbb{R} \rightarrow \mathbb{R}$ be a
piecewise smooth function with discontinuities at $t_1, t_2, \dots$.
Its pointwise derivative $f^{(1)}$ which exists and is continuous
except at $t_1, t_2, \dots$, and its distribution derivative
$f^\prime$ in Schwartz's sense are, as well known, related by
\begin{equation}\label{distrib}
f^\prime (t) = f^{(1)} (t) + \left(f(t_1 +) - f(t_1 -) \right)
\delta (t- t_1) + \left(f(t_2 +) - f(t_2 -) \right) \delta (t- t_2)
+ \dots
\end{equation}
where
\begin{itemize}
\item $f(\tau +) = \lim_{t \downarrow \tau} f(t)$,
$f(\tau -) = \lim_{t \uparrow \tau} f(t)$,
\item $\delta$ is the Dirac delta function.
\end{itemize}
A huge literature\footnote{See the excellent account due to
Basseville and Nikiforov \cite{basseville} for more details.} has
been devoted to the detection of $t_1, t_2, \dots$, which is a
crucial question in signal processing, in diagnosis, and in many
other fields of engineering and applied sciences, where those
discontinuities are often called {\em change-points} or {\em abrupt
changes}.\footnote{The most popular terminology in French is {\em
ruptures}.} Difficulties are stemming from
\begin{itemize}
\item corrupting noises which are blurring the discontinuities,
\item the combined need of
\begin{itemize}
\item fast online calculations,
\item a feasible software and/or hardware implementation.
\end{itemize}
\end{itemize}
Most of the existing literature is based on statistical tools (see,
for instance, \cite{basseville,brodsky1,brodsky2,csorgo} and the
references therein).

The origin of our algebraic viewpoint lies in the references
\cite{esaim,garnier} which are devoted to the parametric
identification of linear systems in automatic
control.\footnote{Change-point detection has also been studied in
\cite{belkoura} via tools stemming from \cite{esaim,garnier}, but in
a quite different manner when compared to us.} We employ elementary
techniques stemming from operational calculus,\footnote{Mikusinski's
foundation \cite{miku1,miku2} of operational calculus, which is not
based on the usual Laplace transform, is a better choice for the
connection with the other algebraic tools. Mikusinski's work, which
is a superb example of {\it algebraic analysis}, is too much
neglected in spite of some advertisements like the nice book by
Yosida \cite{yosida}.} differential algebra and noncommutative
algebra. We are replacing Eq. (\ref{distrib}) by its operational
analogue which is easier to handle. Restricting ourselves to
solutions of operational linear differential equations with rational
coefficients lead to noncommutative rings of linear differential
operators. By representing a change-point by a {\em delay operator},
i.e., an operational exponential, Sect. \ref{prelim} concludes with
the {\em identifiability} of change-points, i.e., the possibility of
expressing them via measured data.\footnote{In the context of
constant linear control systems with delays, which bears some
similarity with what is done here, the identification of delays has
also been tackled in \cite{belkoura2,ollivier,rudolph} via
techniques from \cite{esaim,garnier}.} Higher order change-points,
i.e., discontinuities of derivatives of various orders are briefly
discussed in Sect. \ref{higher}. Sect. \ref{implement} presents
several successful numerical experiments,\footnote{Let us emphasize
that our techniques have already been applied in some concrete
case-studies, where the signals to be processed are stemming from
either biology \cite{spike,zoran} or finance
\cite{finance,finance2}.} which
\begin{itemize}
\item exhibit good robustness properties with respect to several
types of additive and multiplicative corrupting noises;
\item indicate that our approach is still
valid outside of its full mathematical
justification.\protect\footnote{It goes without saying that this
Section, which is mainly descriptive, is not intended to be fully
rigorous.}
\end{itemize}
Preliminary results may be found in \cite{mexico} and
\cite{gretsi,ajaccio}.

\vspace{0.2cm}

\noindent{\bf Acknowledgement}. The authors would like to thank
anonymous referees for several most helpful comments.

\section{Algebra via operational calculus}
\label{prelim}

\subsection{Differential equations} \label{sec:2} Take a
commutative field $k_0$ of characteristic zero. The field $k_0(s)$
of rational functions over $k_0$ in the indeterminate $s$ is
obviously a {\em differential field} with respect to the derivation
$\frac{d}{ds}$ and its subfield of constants is $k_0$ (cf.
\cite{chambert,singer}). Write $k_0(s) [\frac{d}{ds}]$ the
noncommutative ring of linear differential operators of the form
\begin{equation}\label{varpi}
\sum_{\tiny \mbox{\rm finite}} \varrho_\alpha
\frac{d^\alpha}{ds^\alpha} \in k_0 (s) \left[\frac{d}{ds}\right],
\quad \varrho_\alpha \in k_0(s)
\end{equation}
We know that $k_0(s) [\frac{d}{ds}]$ is a left and right principal
ideal ring (cf. \cite{McC,singer}\footnote{Note that \cite{singer}
is not employing, contrarily to \cite{McC}, the usual terminology of
ring and module theory.}). Any signal $x$ is assumed\footnote{See
also \cite{arima,mexico,mboup}.} here to be {\em operationally
holonomic}, i.e., to satisfy a linear differential equation with
coefficients in $k_0(s)$: there exists a linear differential
operator $\varpi \in k_0(s) [\frac{d}{ds}]$, $\varpi \neq 0$, such
that $\varpi x = 0$.
\begin{remark}
Let us explain briefly this assumption. We consider only {\em
holonomic} time functions $z(t)$, i.e., time functions which satisfy
linear differential equations with polynomial coefficients:
$$\left( \sum_{\iota = 0}^{N} p_\iota (t) \frac{d^\iota}{dt^\iota}
\right) z = 0, \quad p_\iota \in \mathbb{C}[t], \quad t\geq 0$$ The
corresponding operational linear differential equation reads (cf.
\cite{yosida})
$$
\left( \sum_{\iota = 0}^{N} p_\iota \left(- \frac{d}{ds}\right)
s^\iota \right) \hat{z} = I(s)
$$
where $ I \in \mathbb{C}[s]$ depends on the initial conditions. A
homogeneous linear differential equation is obtained by
differentiating both sides of the previous equation enough times
with respect to $s$.
\end{remark}
Let $\bar{K}$ be the algebraic closure of $k_0(s)$: $\bar{K}$ is
again a differential field with respect to $\frac{d}{ds}$ and its
subfield of constants is the algebraic closure $\bar{k}_0$ of $k_0$.
It is known that $x$ belongs to a {\em Picard-Vessiot extension} of
$\bar{K}$ (cf. \cite{chambert,singer}).

\begin{remark}
Holonomic functions play an important r\^{o}le in many parts of
mathematics like, for instance, combinatorics (see, e.g.,
\cite{stanley}).
\end{remark}

\subsection{Annihilators}\label{annihilator} Consider now the left
$k_0(s) [\frac{d}{ds}]$-module $\mathcal{M}$ spanned by a finite set
$\{x_\iota | \iota \in I \}$ of such signals. Any $x_\iota$ is a
{\em torsion} element (cf. \cite{McC}) and therefore $\mathcal{M}$
is a {\em torsion module}.\footnote{Such a module is called a {\em
differential module} in \cite{singer}.} The {\em annihilator} $A_I$
of $\{x_\iota | \iota \in I \}$ is the set of linear differential
operators $\varpi_I \in k_0(s) [\frac{d}{ds}]$ such that, $\forall ~
\iota \in I$, $\varpi_I x_\iota = 0$. It is a left ideal of $k_0(s)
[\frac{d}{ds}]$ and it is therefore generated by a single element
$\Delta \in k(s) [\frac{d}{ds}]$, $\Delta \neq 0$, which is called a
{\em minimal} annihilator of $\{x_\iota | \iota \in I \}$. It is
obvious that $\Delta$ is annihilating any element belonging to the
$k_0$-vector space $\mbox{\rm span}_{k_0} (x_\iota | \iota \in I)$.
The next result is straightforward:
\begin{lemma}\label{unique}
Let $\Delta_1$, $\Delta_2$, $\Delta_1 \neq \Delta_2$, be two minimal
annihilators. There exists $\rho \in k_0 (s)$, $\rho \neq 0$, such
that $\Delta_1 = \rho \Delta_2$.
\end{lemma}
We will say that the minimal annihilator is unique up to left
multiplications by nonzero rational functions.

A rational function $\frac{p}{q}$, $p, q \in k_0[s]$, $q \neq 0$, is
said to be {\em proper} (resp. {\em strictly proper}) if, and only
if, the $\mbox{\rm d}^\circ p \leq \mbox{\rm d}^\circ q$ (resp.
$\mbox{\rm d}^\circ p < \mbox{\rm d}^\circ q$). A differential
operator (\ref{varpi}) is said to be {\em proper} (resp. {\em
strictly proper}) if, and only if, any $\varrho_\alpha$ is proper
(resp. strictly proper). The next result is an obvious corollary of
Lemma \ref{unique}:
\begin{corollary}\label{proper}
It is possible to choose an annihilator, which is minimal or not, in
such a way that it is proper (resp. strictly proper).
\end{corollary}

A rational function $\frac{p}{q}$, $p, q \in k_0[s]$, $q \neq 0$, is
said to be in a {\em finite integral form} (resp. {\em strictly
finite integral form}) if, and only if, it belongs to $k_0
[\frac{1}{s}]$ (resp. $\frac{1}{s} k [\frac{1}{s}]$). A differential
operator (\ref{varpi}) is said to be in a {\em finite integral form}
(resp. {\em strictly finite integral form}) if, and only if, any
$\varrho_\alpha$ is in a finite integral form (resp. strictly finite
integral form). Consider a common multiple $m \in k[s]$ of the
denominators of the $\varrho_\alpha$'s. The operator $s^{- N} m
\varpi$ is in a (strictly) finite integral form for large enough
values of the integer $N \geq 0$.

\begin{corollary}\label{int}
It is possible to choose an annihilator, which is minimal or not, in
such a way that it is in a finite integral form (resp. strictly
finite integral form).
\end{corollary}

\subsection{Delay operators} Let $k/k_0$ be a transcendental field
extension. The field $k(s)$ of rational functions over $k$ in the
indeterminate $s$ is again a differential field with respect to
$\frac{d}{ds}$ and its subfield of constants is $k$. The
noncommutative ring $k (s) [\frac{d}{ds}]$ of linear differential
operators is defined as in Sect. \ref{sec:2}. Pick up an element
$t_r \in k$, called {\em delay}, which is transcendental over $k_0$.
Write the {\em delay operator} with its classic exponential notation
$e^{- t_r s}$ (cf. \cite{pol}), as it satisfies the differential
equation $\left(\frac{d}{ds} + t_r \right) e^{- t_r s} = 0$.
According to Sect. \ref{annihilator}, the differential operator
$\frac{d}{ds} + t_r \in k (s) [\frac{d}{ds}]$ is a minimal
annihilator of $e^{- t_r s}$.

\subsection{Identifiability of the delay}\label{identifiability}
\subsubsection{Main result} Let $\varpi_1, \varpi_2
\in k_0(s) [\frac{d}{ds}]$ be minimal annihilators of two signals
$x_1$, $x_2$, $x_1 x_2 \neq 0$. Introduce the quantity
\begin{equation}\label{obs}
X = x_1 + x_2 e^{- t_r s}
\end{equation}
Multiplying on the left both sides of Eq. (\ref{obs}) by $\varpi_1$
yields $\varpi_1 X =  \varpi_1  x_2 e^{- t_r s}$. Thus
$$
\varpi_1 X e^{t_r s} =  \varpi_1  x_2
$$
and
\begin{equation}\label{exp}
\varpi_2^\prime \varpi_1 X e^{t_r s} = 0
\end{equation}
where $\varpi_2^\prime \in k_0 (s) [\frac{d}{ds}]$ is a minimal
annihilator of $\varpi_1  x_2$. The next proposition follows at
once:
\begin{proposition}\label{main}
Eq. \eqref{exp} is equivalent to
\begin{equation*}\label{alg}
\sum_{\tiny{\mbox{\rm finite}}} t^{\nu}_{r} \left(\pi_\nu X \right)
= 0, \quad \pi_\nu \in k_0(s) \left[\frac{d}{ds}\right]
\end{equation*}
where at least one $\pi_\nu$ is not equal to $0$.
\end{proposition}
Write $k_0 (s) \langle X \rangle$ the differential overfield of $k_0
(s)$ generated by $X$.
\begin{corollary}\label{alg1}
$t_r$ in Eq. \eqref{obs} is algebraic over the differential field
$k_0 (s) \langle X \rangle$.
\end{corollary}
\begin{remark}\label{remark}
Assume that the quantity $X$ is {\em measured}, i.e., there exists a
{\em sensor} which gives at each time instant its numerical value in
the time domain. Then, according to the terminology in \cite{easy},
Corollary \ref{alg1} may be rephrased by saying that $t_r$ is {\em
algebraically identifiable}.
\end{remark}

\subsubsection{First example}\label{first}
Set $x_i = \sum_{\nu_i = 0}^{N_i} \frac{\gamma_{\nu_i}}{s^{\nu_i +
1}}$, $i = 1, 2$, $\gamma_{\nu_i} \in k_0$, where $N_i$ is a known
non-negative integer.\footnote{The coefficients $\gamma_{\nu_i}$ are
not necessarily known.} Then $\frac{d^{N_i + 2}}{ds^{N_i + 2}}
s^{N_i + 1}$ is a minimal annihilator of $x_i$. It follows at once
that Proposition \ref{main} and Corollary \ref{alg1} apply to this
case.

Straightforward calculations demonstrate that $t_r$ is the unique
solution of an equation of the form
\begin{equation}\label{mult}
\mathfrak{p}  \left( \frac{d}{ds} + t_r \right)^\varrho X = 0
\end{equation}
where $\mathfrak{p} \in k_0(s) [\frac{d}{ds}]$, $1 \leq \varrho \leq
N_2$.

\subsubsection{Second example}\label{second}
Assume that $x_2 = \frac{\bar{\gamma}_{2}}{s^{N_2 + 1}}$,
$\bar{\gamma}_{2} \in k_0$, in Sect. \ref{first}. Multiply both
sides of Eq. \eqref{obs} by $(\frac{d}{ds} + t_r) s^{N_2 + 1}$
yields
$$
(\frac{d}{ds} + t_r) s^{N_2 + 1} X = (\frac{d}{ds} + t_r) s^{N_2 +
1} x_1
$$
Eq. \eqref{mult} becomes
\begin{equation}\label{mono}
\pi_1 (\frac{d}{ds} + t_r) s^{N_2 + 1} X =  0
\end{equation}
where $\pi_1 \in k_0 (s) [\frac{d}{ds}]$ is a minimal annihilator of
$(\frac{d}{ds} + t_r) s^{N_2 + 1} x_1$.
\begin{proposition}\label{lin1}
$t_r$ satisfies an algebraic Equation \eqref{mono} of degree $1$.
\end{proposition}
\begin{remark}\label{momo}
If $X$ is measured as in Remark \ref{remark}, then, according to the
terminology in \cite{esaim,garnier}, Proposition \ref{lin1} may be
rephrased by saying that $t_r$ is {\em linearly identifiable}.
\end{remark}

\subsubsection{Third example}\label{rat}
Assume in Eq. \eqref{obs} that $x_2 = \frac{a}{b} \in k_0 (s)$, $a,
b \in k_0 [s]$, $(a, b) = 1$, is a known rational function, i.e.,
\begin{equation}\label{const}
X = x_1 +  \frac{a}{b} e^{- t_r s}
\end{equation}
Multiplying both sides by $(\frac{d}{ds} + t_r) \frac{b}{a}$ yields
$(\frac{d}{ds} + t_r) \frac{b}{a} X = (\frac{d}{ds} + t_r)
\frac{b}{a} x_1$. Since $t_r$ is constant, there exists an
annihilator $\pi \in k_0 (s) [\frac{d}{ds}]$ of $(\frac{d}{ds} +
t_r) \frac{b}{a} x_1$, i.e.,
\begin{equation}\label{const1}
\left(\pi \frac{b}{a} X\right) t_r + \pi \frac{d}{ds}
\left(\frac{b}{a} X \right) = 0
\end{equation}
\begin{proposition}\label{lin}
$t_r$ in Eq. \eqref{const} satisfies an algebraic Equation
\eqref{const1} of degree $1$.
\end{proposition}
\begin{remark}
If $X$ is measured as in Remark \ref{remark}, then, according to
Remark \ref{momo}, Proposition \ref{lin} may be rephrased by saying
that $t_r$ is linearly identifiable.
\end{remark}

\section{Higher order change-points}\label{higher}
Take again as in the introduction a piecewise smooth function $f$,
which is now assumed to be $C^n$, $n \geq 0$, i.e., $f$ and its
pointwise derivatives up to order $n$ are continuous. We might be
interested in the discontinuities of its $(n+1)^{th}$ order
pointwise derivative, which are called {\em change-points}, or {\em
abrupt changes}, {\em of order $n + 1$}.

By replacing Eq. \eqref{obs} by
\begin{equation*}\label{obser}
s^{(n+1)} X = x_1 + x_2 e^{- t_r s}
\end{equation*}
it is straightforward to extend all the results of Sect.
\ref{identifiability} to higher order change-points.

\section{Some numerical experiments}\label{implement}
\subsection{General principles}\label{justi}
From now on $k_0$ is a subfield of $\mathbb{R}$, $\mathbb{Q}$ for
instance. We utilize the calculations of Sect. \ref{second} like
follows:
\begin{itemize}
\item Multiplying both sides of Eq. \eqref{mono} by $s^{- N}$, where
$N > 0$ is large enough, yields
\begin{equation}\label{iterint}
s^{- N} \pi_1 \left(\frac{d}{ds} + t_r\right) s^{N_2 + 1} X =  0
\end{equation}
where $s^{- N} \pi_1 (\frac{d}{ds} + t_r) s^{N_2 + 1}$ is a strictly
integral operator.
\item Going back to the time domain is achieved via the classic
rules of operational calculus \cite{miku1,miku2,yosida}, where
$\frac{d^\nu}{ds^\nu}$ corresponds to the multiplication by $(-
t)^\nu$.
\item $x_1 = \sum_{\nu_1 = 0}^{N_1} \frac{\gamma_{\nu_1}}{s^{\nu_1 +
1}}$ and $x_2 = \frac{\bar{\gamma}_{2}}{s^{N_2 + 1}}$ correspond in
the time domain to the polynomial functions $\sum_{\nu_1 = 0}^{N_1}
\frac{\gamma_{\nu_1} t^{\nu_1}}{\nu_1 !}$ and
$\frac{\bar{\gamma}_{2} t^{N_2}}{N_2 !}$.
\item Those time functions are assumed to approximate on a ``short''
time interval the signal where change-points have to be detected.
\item Consider the numerical value $v$ taken by the time analogue of the left side
of Eq. \eqref{iterint} when the value given to $t_r$ is the middle
of a given ``short'' time window. If $v$ is ``close'' to $0$, then
we say that the middle of the time window is a change-point.
\item This time window is sliding in order to capture the various
change-points, which are assumed to be not too ``close'', i.e., the
distance between two consecutive change-points is larger than the
time window.
\item The corrupting noises are attenuated by the iterated time integrals
which corresponds in the time domain to the negative power of $s$ in
the left side of Eq. \eqref{iterint}.\footnote{Noises in \cite{ans}
are viewed, via {\em nonstandard analysis}, as quickly fluctuating
phenomena (see also \cite{lobry} for an introductory presentation).
The noises are attenuated by the iterated time integrals, which are
simple examples of low-pass filters (we may also choose, according
to Lemma \ref{int}, more involved low-pass filters (see, e.g.,
\cite{chen})). No statistical tools are required and we are by no
means restricted to Gaussian white noises, like too often in the
engineering studies. Moreover the corrupting noises need not to be
additive. They might also be multiplicative.}
\end{itemize}

\subsection{Examples\protect\footnote{Interested
readers may ask C. Join for the corresponding computer programs
({\tt Cedric.Join}@{\tt cran.uhp-nancy.fr}).}}

\begin{figure}[h]
 \centering
 {\subfigure[\footnotesize Noise-free signal (- -), signal (--)]{\rotatebox{-90}{\includegraphics[width=4cm]{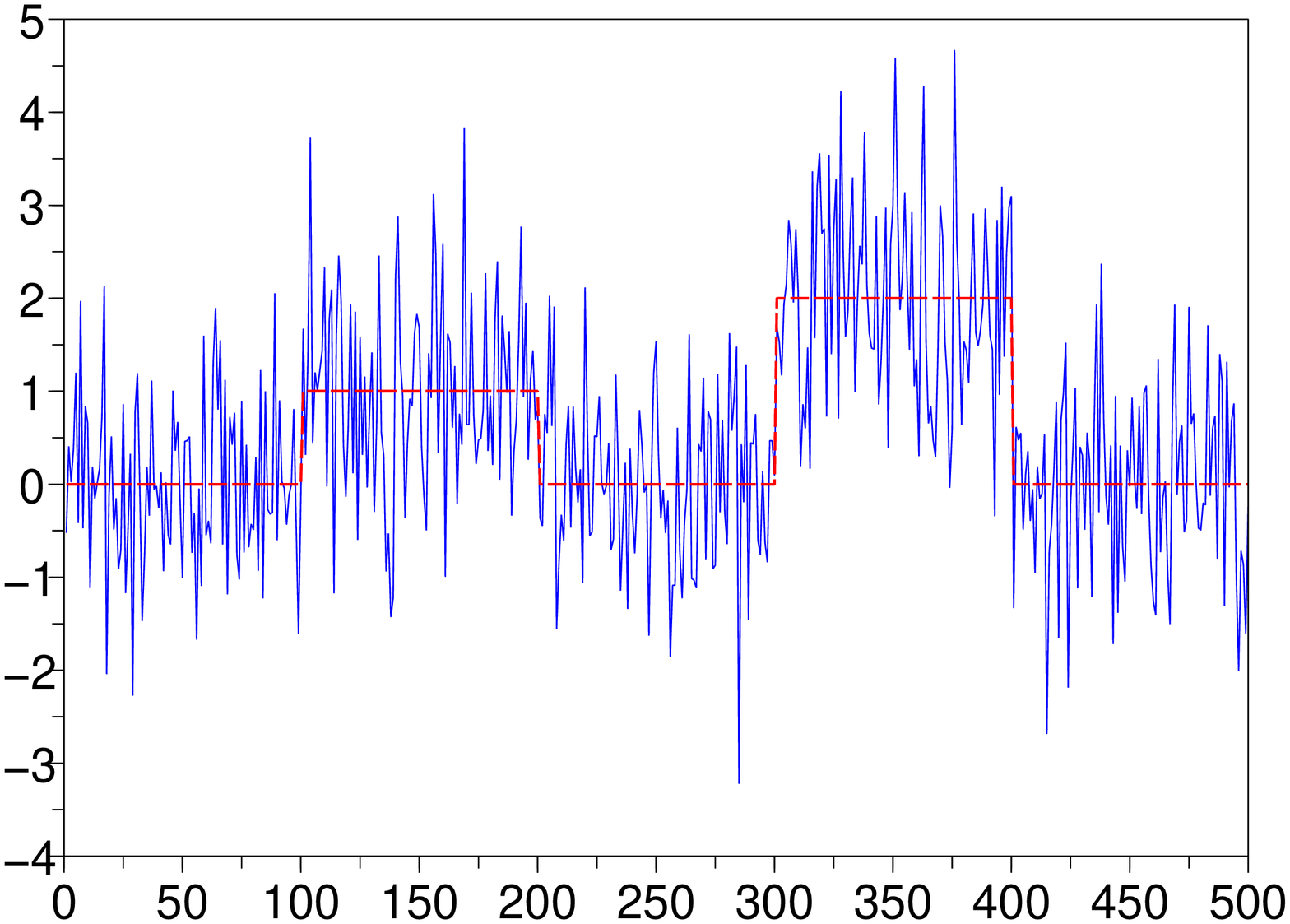}}}}
 {\subfigure[\footnotesize Change-point detection -- Exact (+)]{\rotatebox{-90}{\includegraphics[width=4cm]{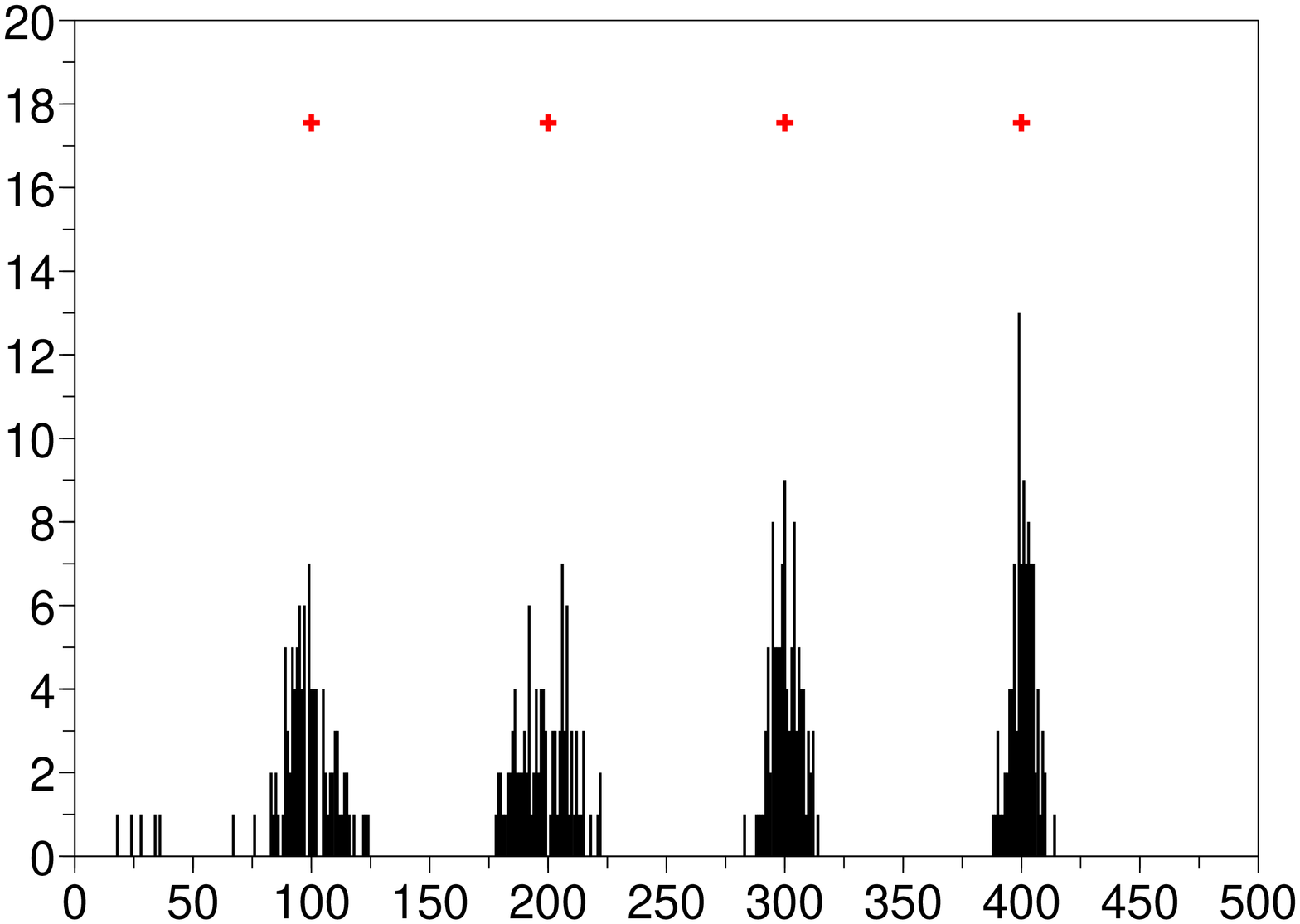}}}}
\caption{Piecewise constant signal -- Normal additive noise -- SNR:
0 db} \label{M1}
\end{figure}
\begin{figure}[h]
 \centering
 {\subfigure[\footnotesize Noise-free signal (- -), signal (--)]{\rotatebox{-90}{\includegraphics[width=4cm]{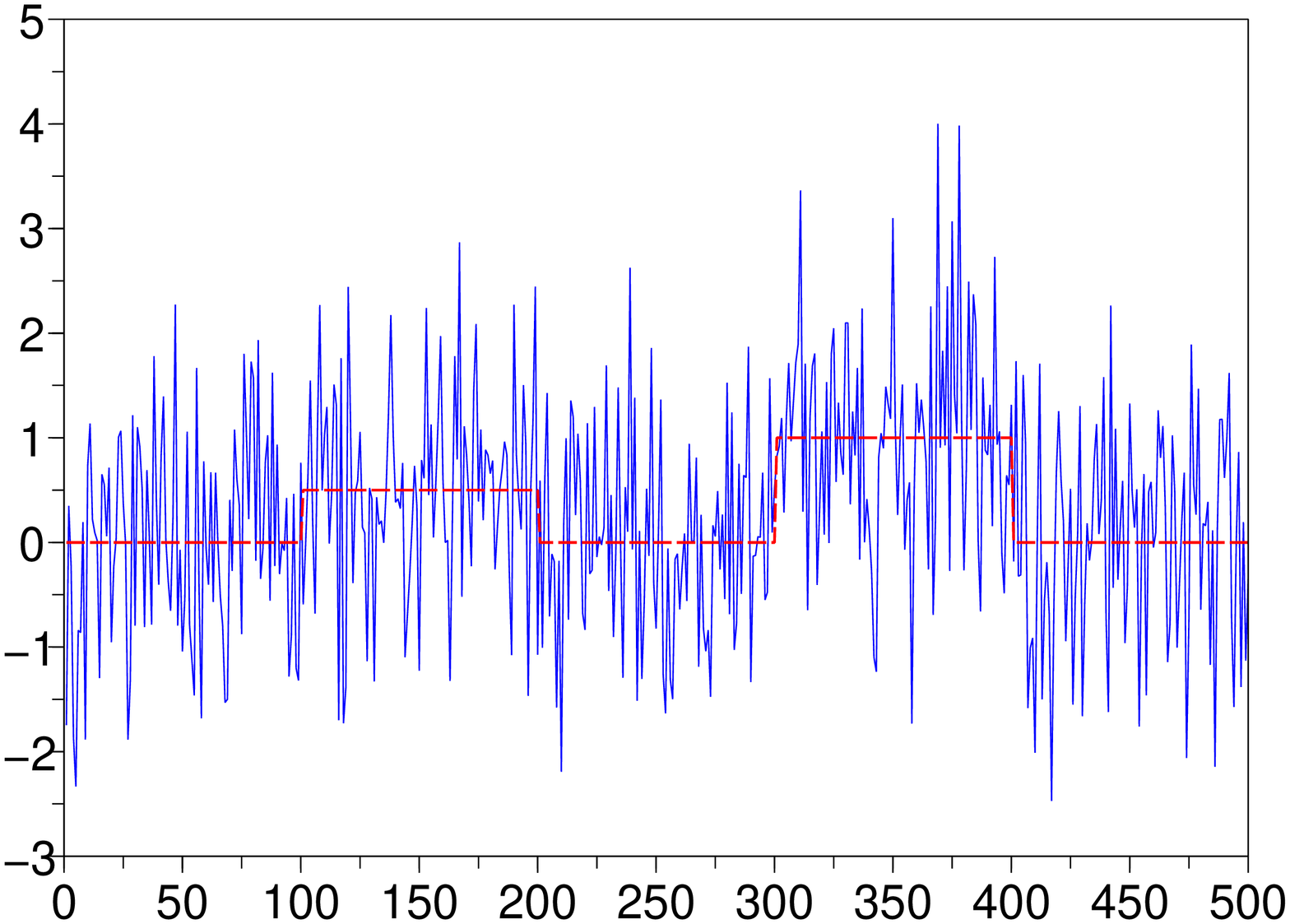}}}}
 {\subfigure[\footnotesize Change-point detection -- Exact (+)]{\rotatebox{-90}{\includegraphics[width=4cm]{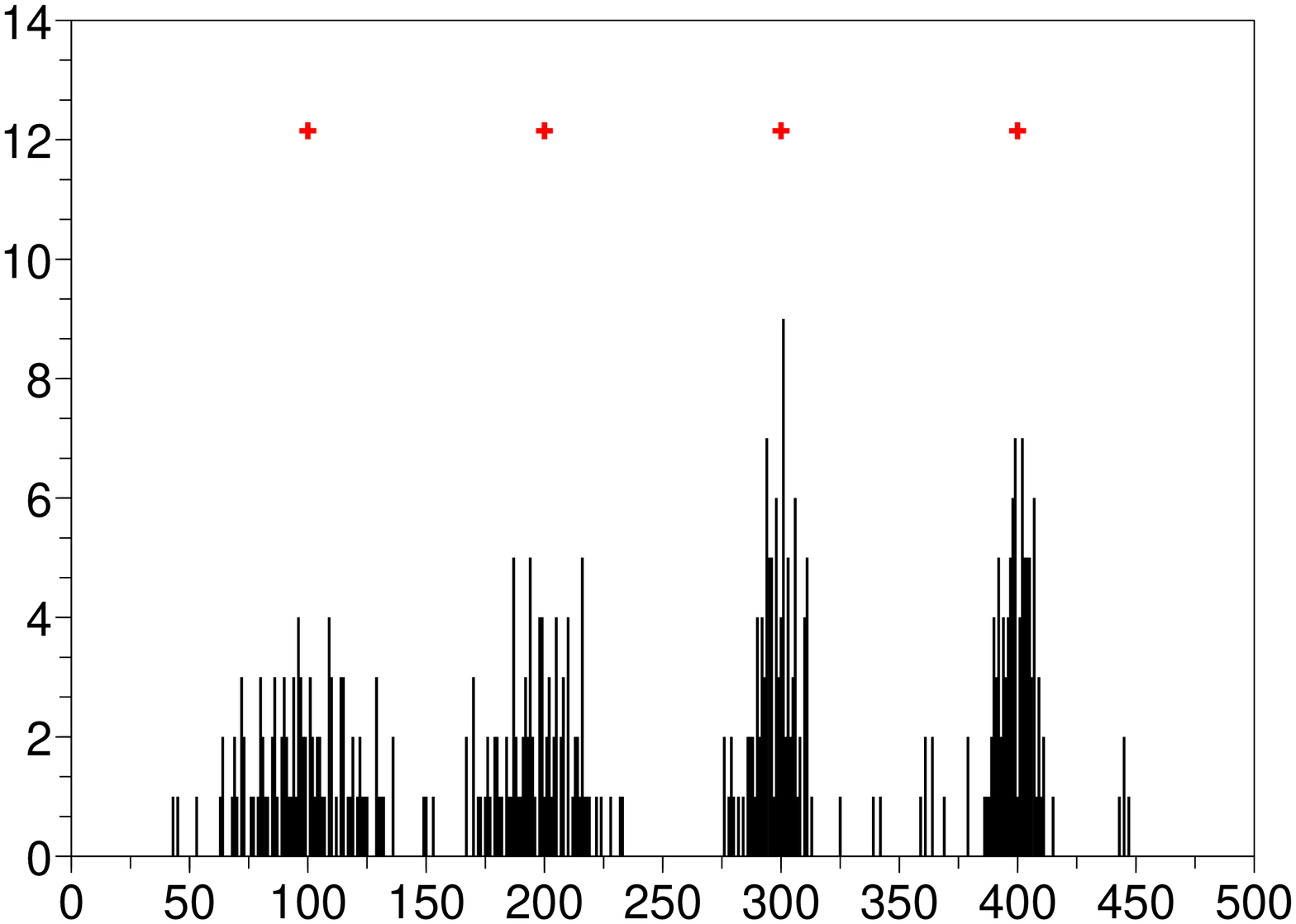}}}}
\caption{Piecewise constant signal -- Normal additive noise -- SNR:
-6 db} \label{M2}
\end{figure}

The following academic examples are investigated:
\begin{itemize}
\item piecewise constant and polynomial real-valued functions,
\item a real-valued sinusoid plus a piecewise constant real-valued
function.
\end{itemize}
The robustness with respect to corrupting noises, which is reported
in Table 1, is tested thanks to several noises, of various
powers,\footnote{We are utilizing the notion of {\em signal-to-noise
ratio}, or {\em SNR}, which is familiar in signal processing (see
{\it Wikipedia}, for instance).} which are of the following types:
\begin{enumerate}
\item additive, zero mean, and either normal, uniform or Perlin,\footnote{{\em
Perlin's noises} \cite{perlin} are quite popular in computer
graphics.}
\item multiplicative, of mean $1$, and uniform.
\end{enumerate}
We finally note that
\begin{itemize}
\item piecewise polynomial functions were difficult to analyze even
via recent techniques like wavelets (see, e.g., \cite{drag,mallat});
\item we do not need any {\it a priori}
knowledge of the upper bound of the number of change-points (see,
e.g., \cite{lavielle,sinica});
\item we are not limited to a given type of noises and we are able to handle
multiplicative noises as well (see, e.g.,
\cite{gij,lebarbier,tourneret});
\item the results remain satisfactory even with a very high noise level
(see Figures \ref{M1} and \ref{M2}).
\end{itemize}

\begin{remark}
The so-called Perlin noises, which are not familiar in signal
processing and in automatic control, contain components which are
obviously not quickly fluctuating. It is all the more remarkable
that our computer simulations are still good, in spite of the fact
that this example goes beyond the theoretical justifications
provided in Sect. \ref{justi}.
\end{remark}


\begin{figure}[h]
 \centering
 {\subfigure[\footnotesize Noise-free signal (- -), signal (--)]{\rotatebox{-90}{\includegraphics[width=4cm]{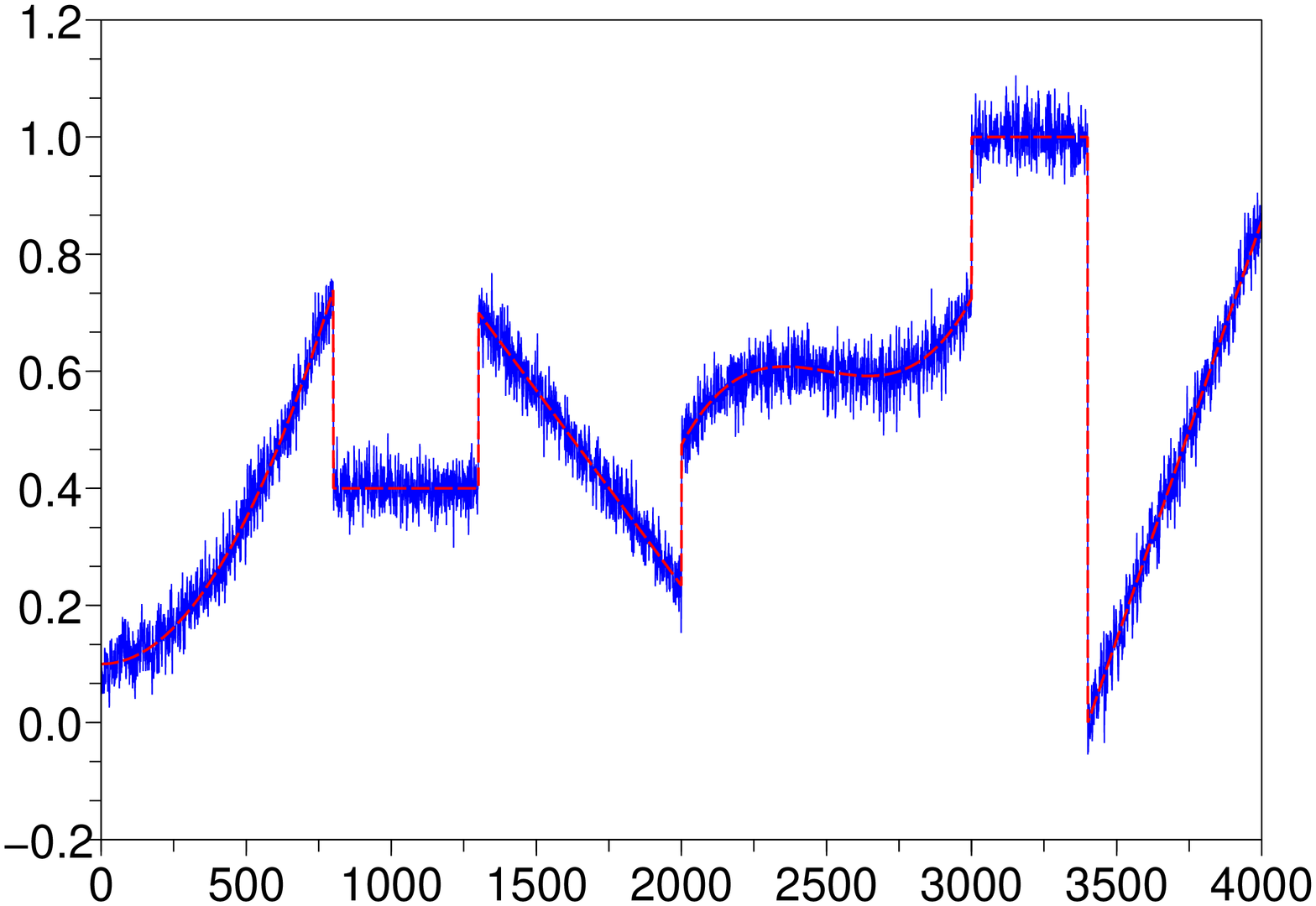}}}}
 {\subfigure[\footnotesize Change-point detection -- Exact (+)]{\rotatebox{-90}{\includegraphics[width=4cm]{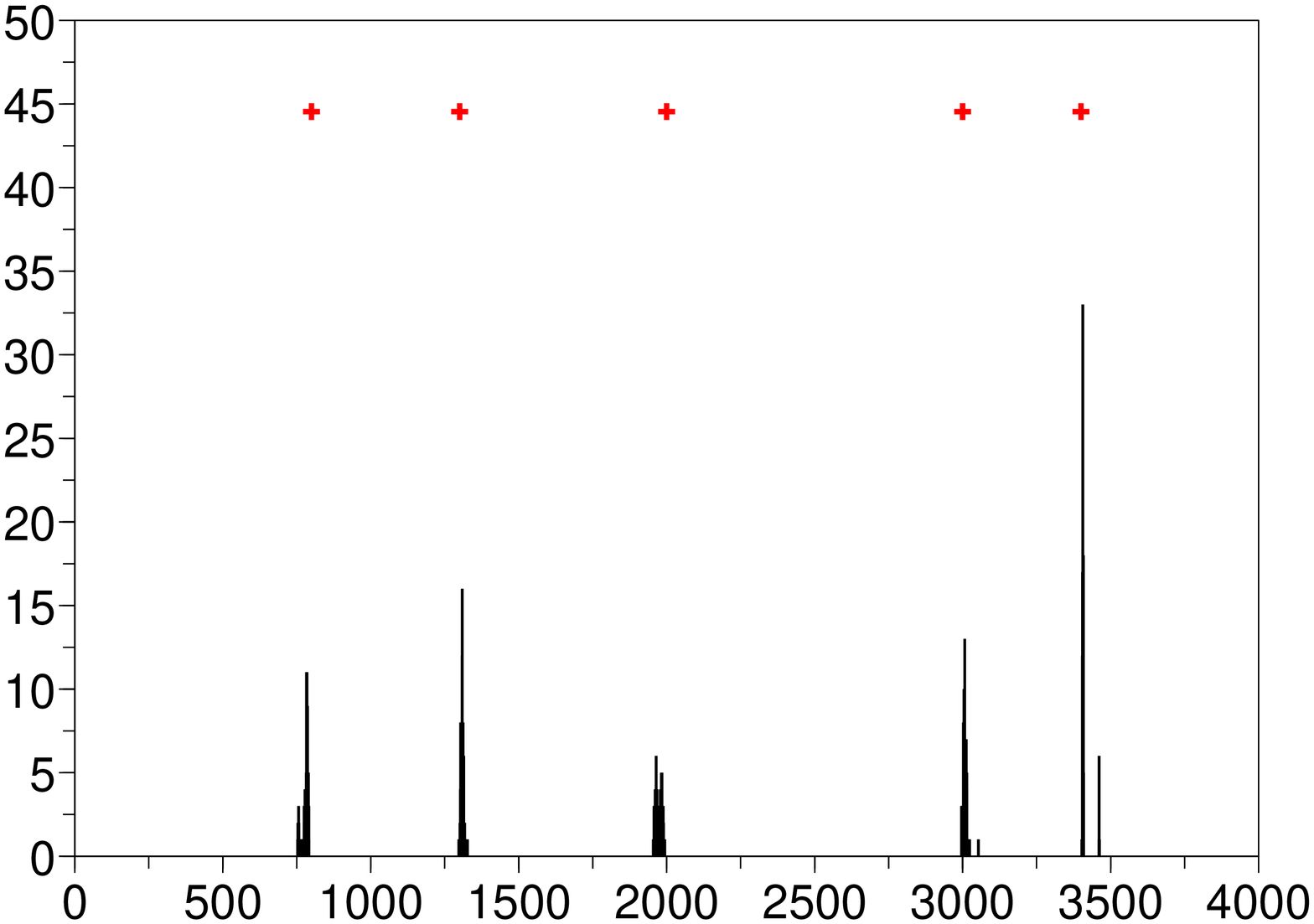}}}}
\caption{Piecewise polynomial signal -- Normal Additive noise --
SNR: 25 db} \label{Bg25}
\end{figure}

\begin{figure}[h]
 \centering
 {\subfigure[\footnotesize Noise-free signal (- -), signal (--)]{\rotatebox{-90}{\includegraphics[width=4cm]{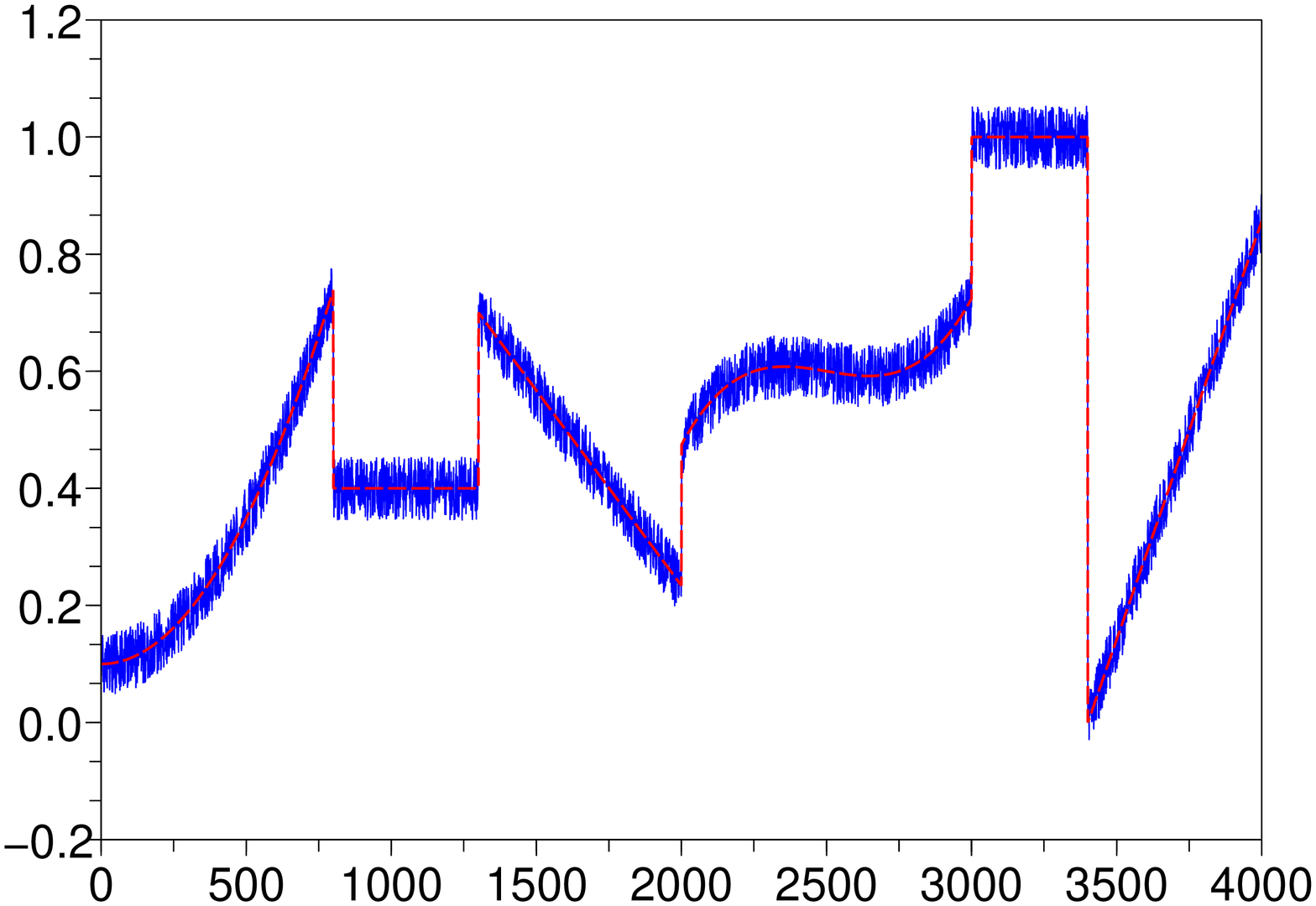}}}}
 {\subfigure[\footnotesize Change-point detection -- Exact (+)]{\rotatebox{-90}{\includegraphics[width=4cm]{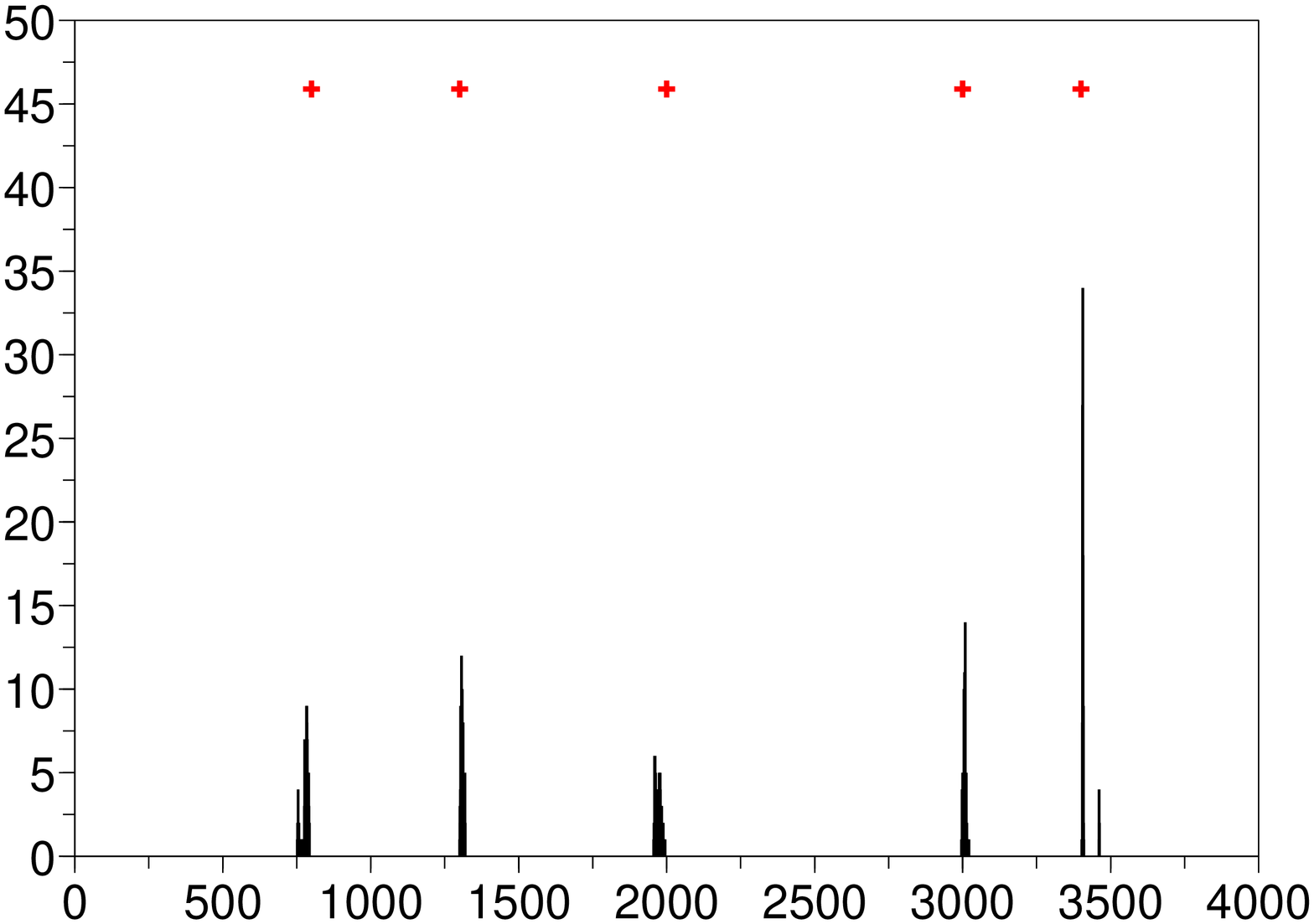}}}}
\caption{Piecewise polynomial signal -- Uniform additive noise  --
SNR: 25 db } \label{Bu25}
\end{figure}

\begin{figure}[h]
 \centering
 {\subfigure[\footnotesize Noisy free signal (- -), signal (--)]{\rotatebox{-90}{\includegraphics[width=4cm]{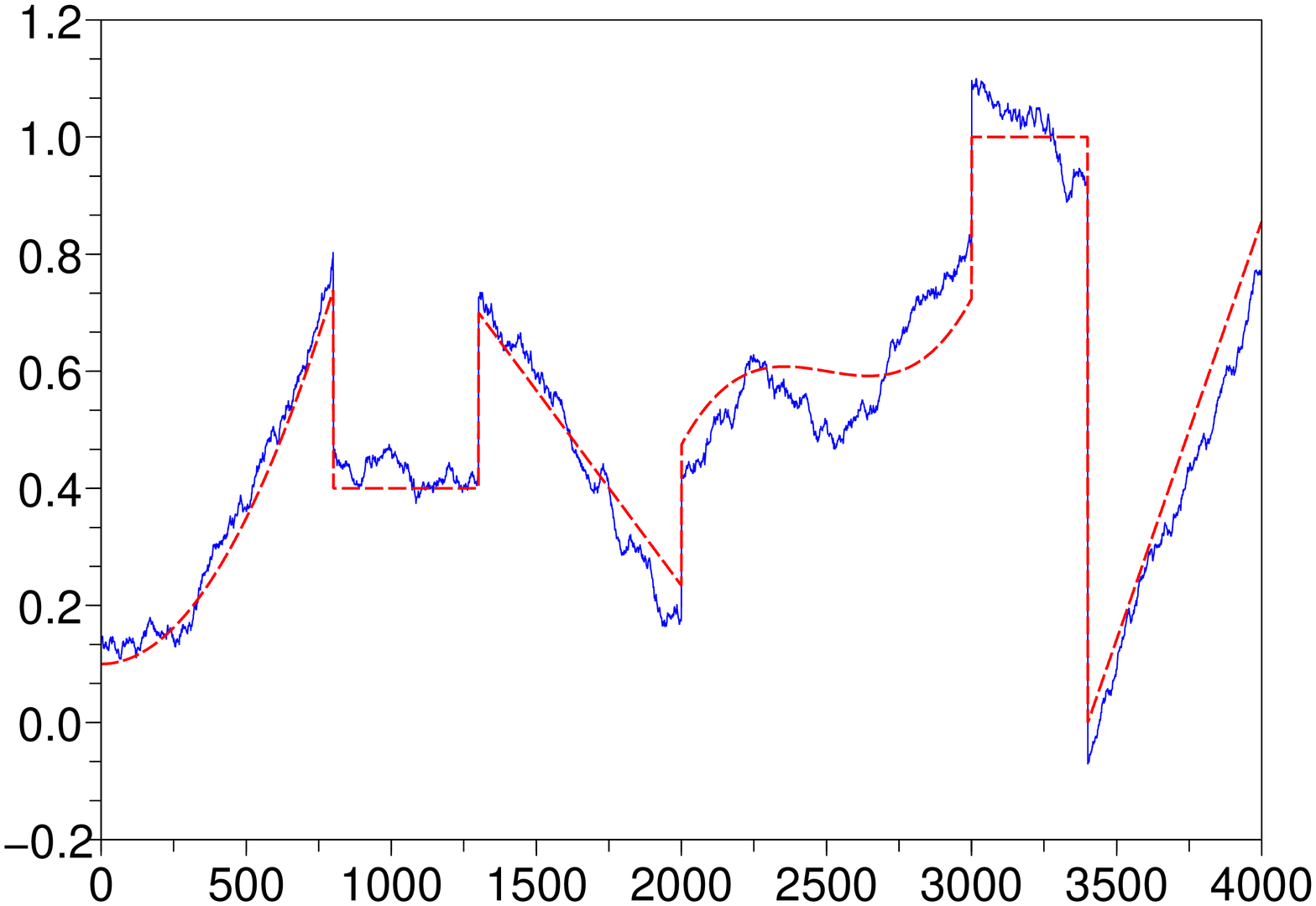}}}}
 {\subfigure[\footnotesize Change-point detection -- Exact (+)]{\rotatebox{-90}{\includegraphics[width=4cm]{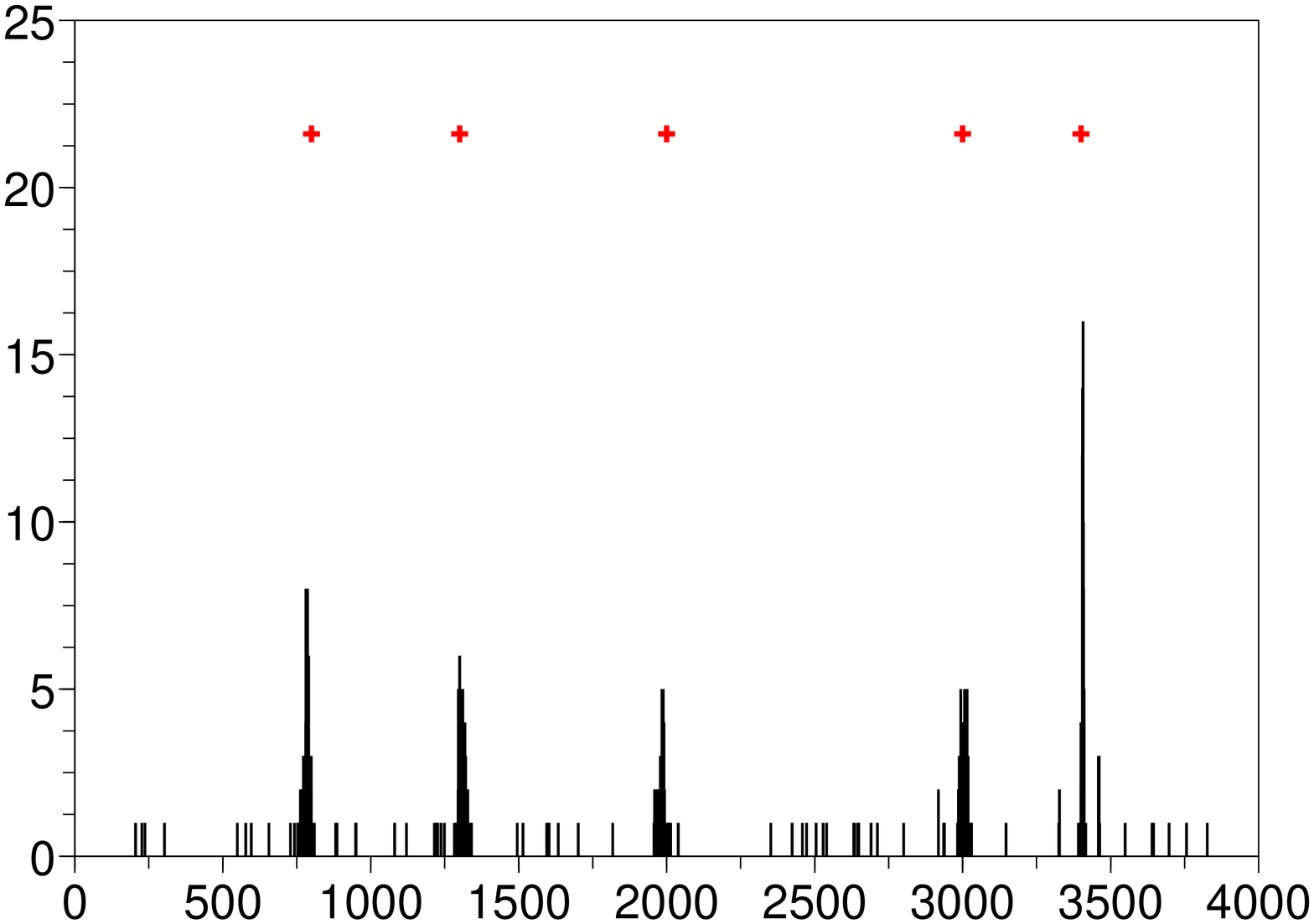}}}}
\caption{Piecewise polynomial signal -- Additive Perlin noise --
SNR: 20 db} \label{Bp20}
\end{figure}

\begin{figure}[h]
 \centering
 {\subfigure[\footnotesize Noise-free signal (- -), signal (--)]{\rotatebox{-90}{\includegraphics[width=4cm]{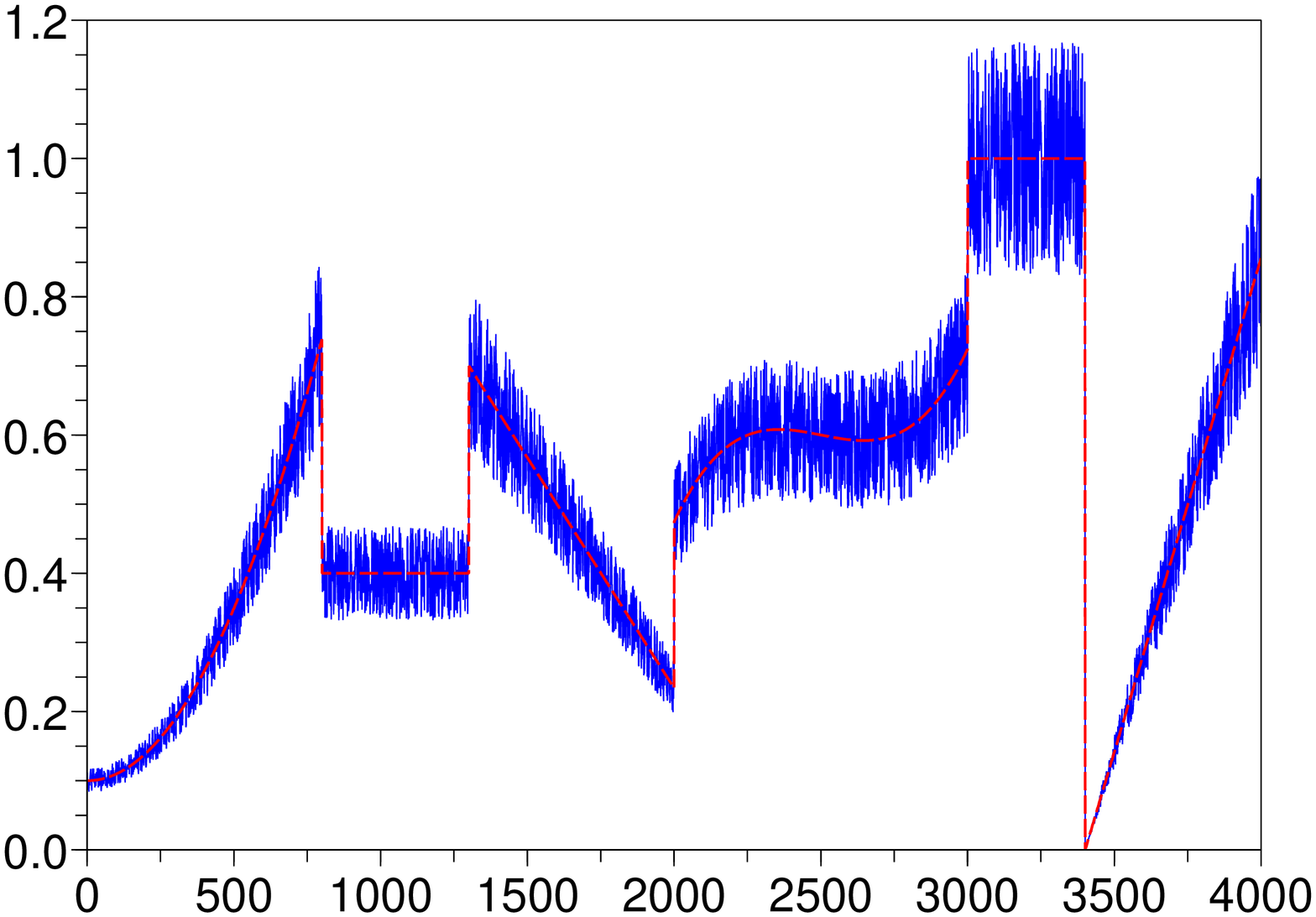}}}}
 {\subfigure[\footnotesize Distribution of change-point detection, exact places (+)]{\rotatebox{-90}{\includegraphics[width=4cm]{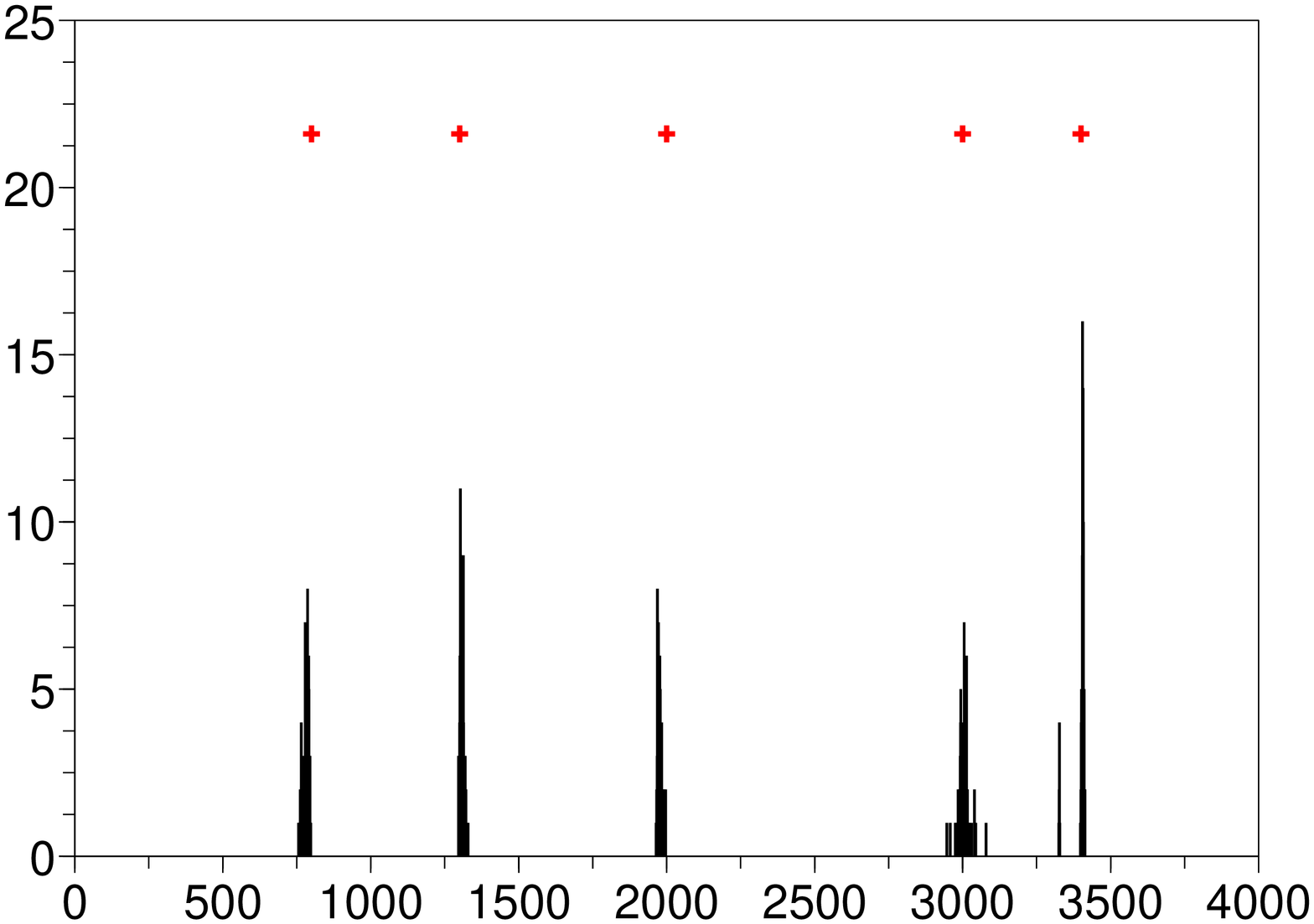}}}}
\caption{Piecewise polynomial signal -- Uniform multiplicative noise
-- SNR: 20 db} \label{Bmu20}
\end{figure}

\begin{figure}[h]
 \centering
 {\subfigure[\footnotesize Noise-free signal (- -), signal (--)]{\rotatebox{-90}{\includegraphics[width=4cm]{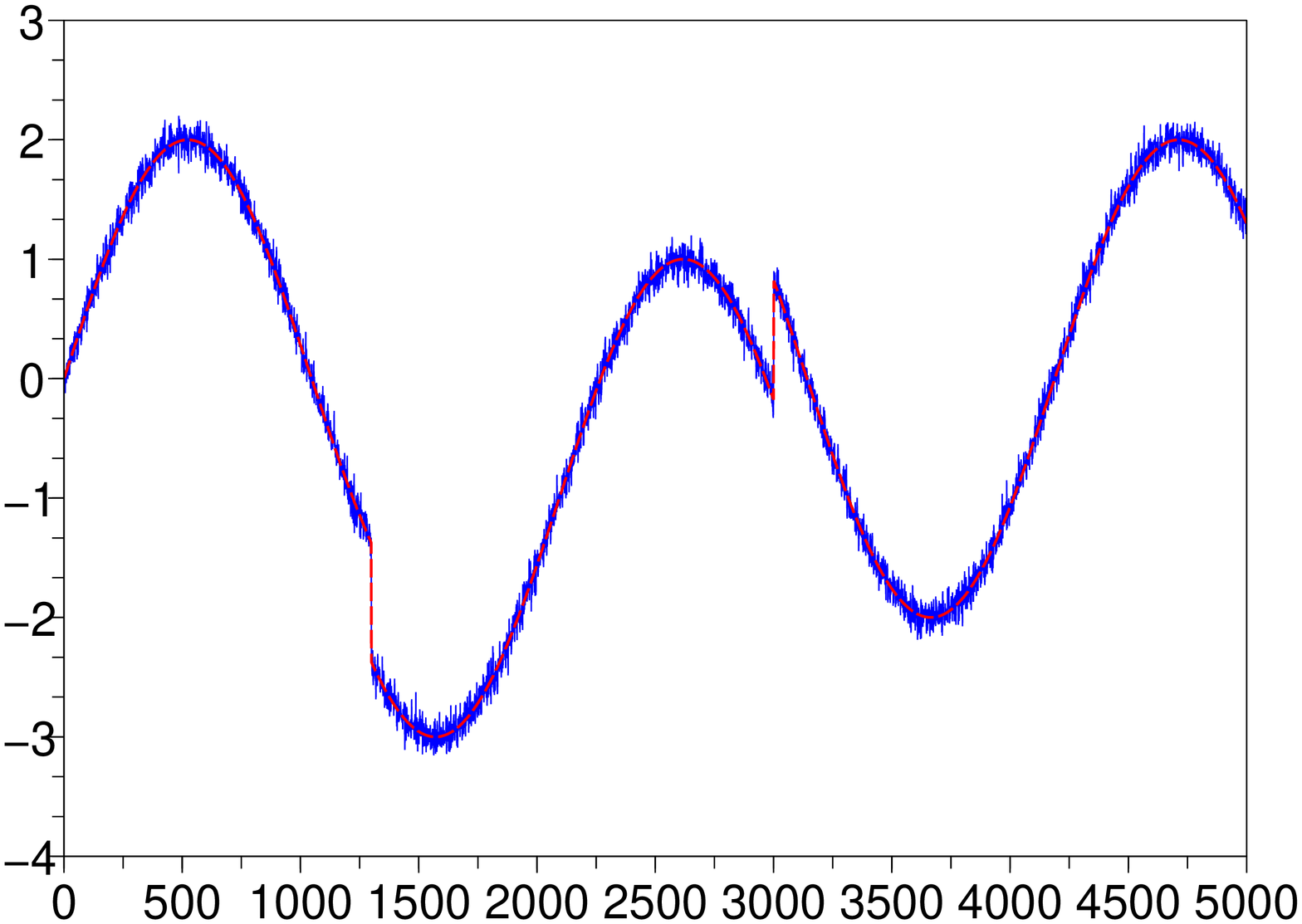}}}}
 {\subfigure[\footnotesize Change-point detection -- Exact (+)]{\rotatebox{-90}{\includegraphics[width=4cm]{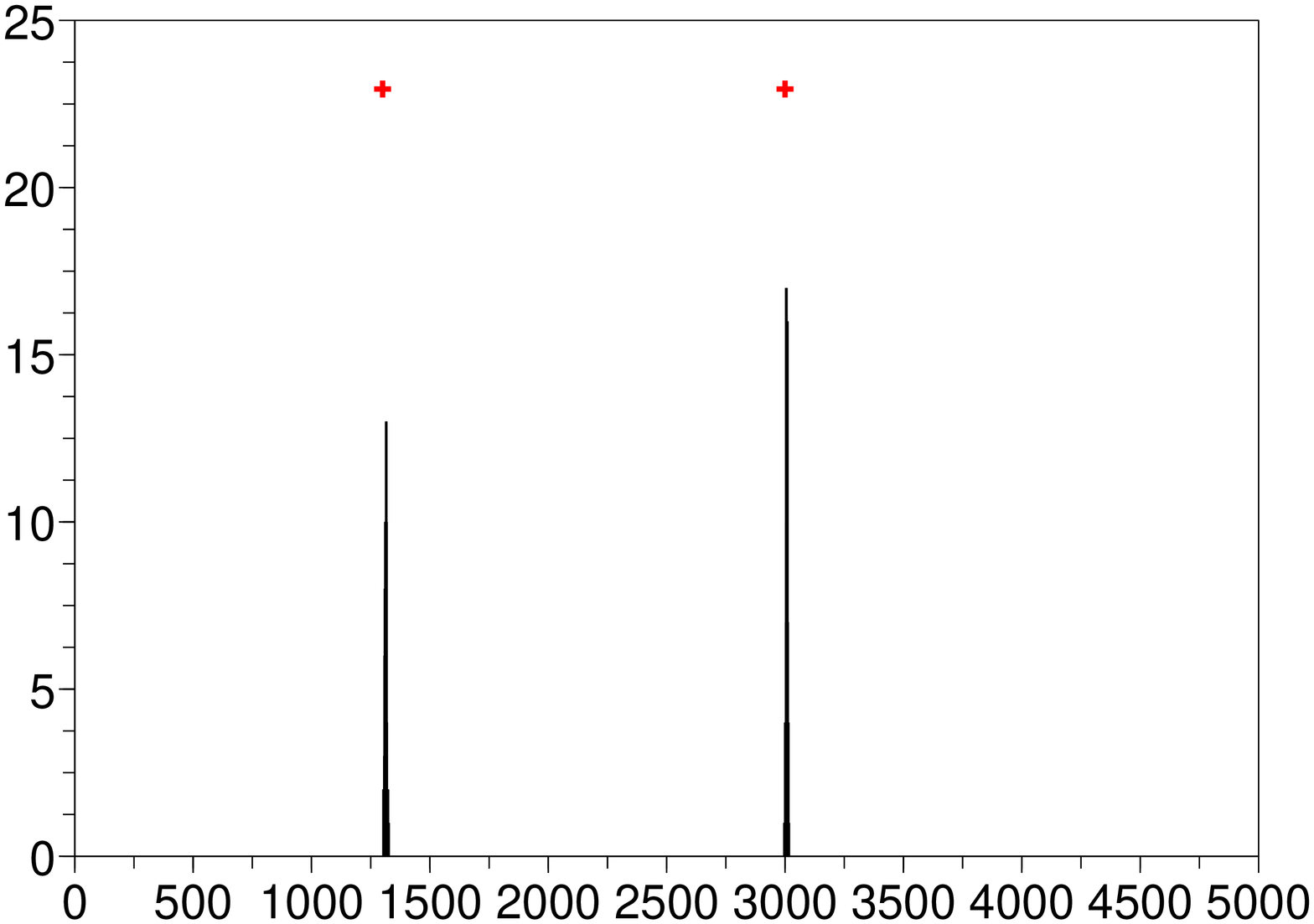}}}}
\caption{Sinusoidal signal -- Normal additive noise -- SNR: 25 db }
\label{Sg25}
\end{figure}

\begin{figure}[h]
 \centering
 {\subfigure[\footnotesize Noise-free signal (- -), signal (--)]{\rotatebox{-90}{\includegraphics[width=4cm]{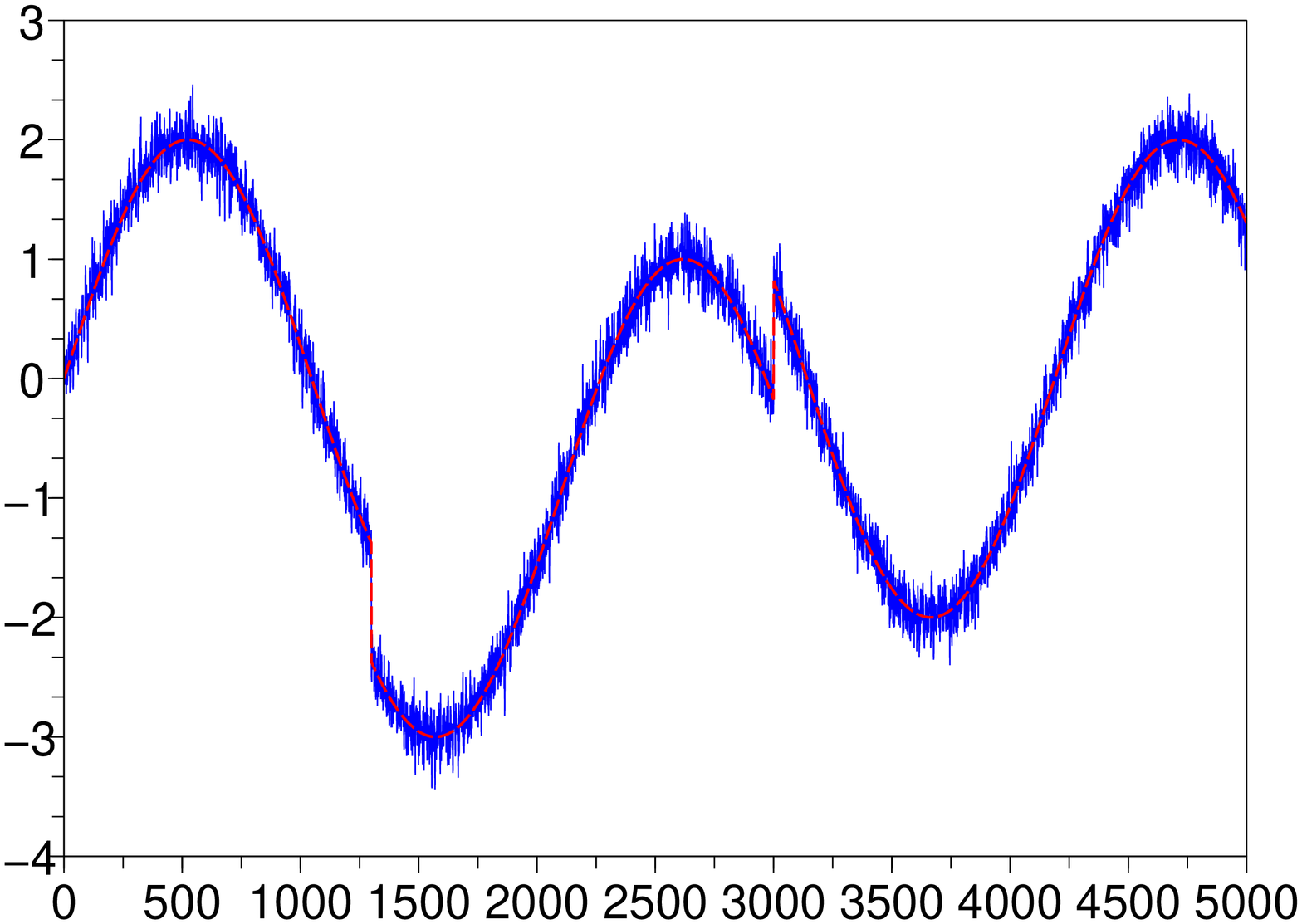}}}}
 {\subfigure[\footnotesize Change-point detection -- Exact (+)]{\rotatebox{-90}{\includegraphics[width=4cm]{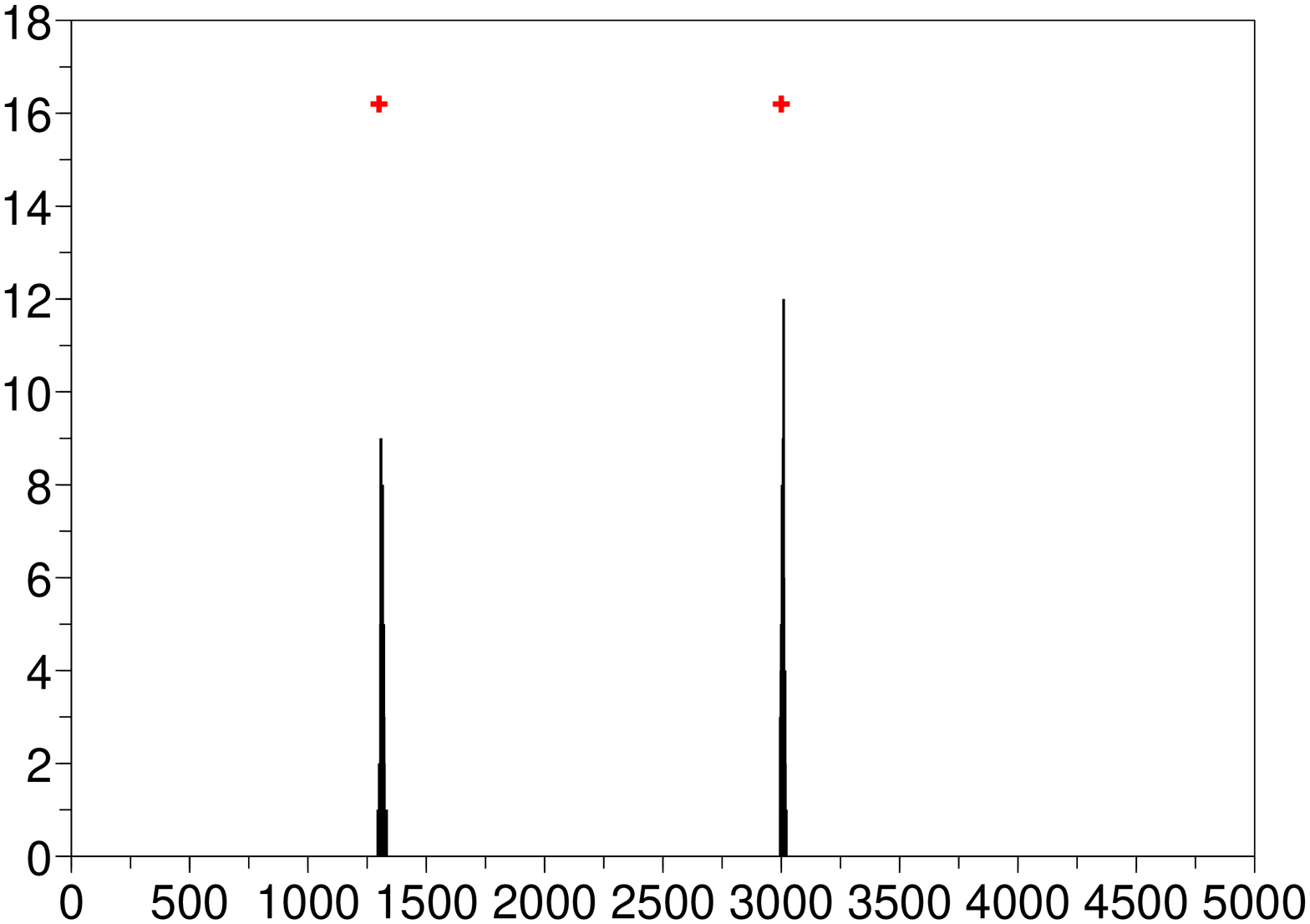}}}}
\caption{Sinusoidal signal -- Normal additive noise -- SNR: 20 db }
\label{Sg20}
\end{figure}

\begin{figure}[h]
 \centering
 {\subfigure[\footnotesize Noise-free signal (- -), signal (--)]{\rotatebox{-90}{\includegraphics[width=4cm]{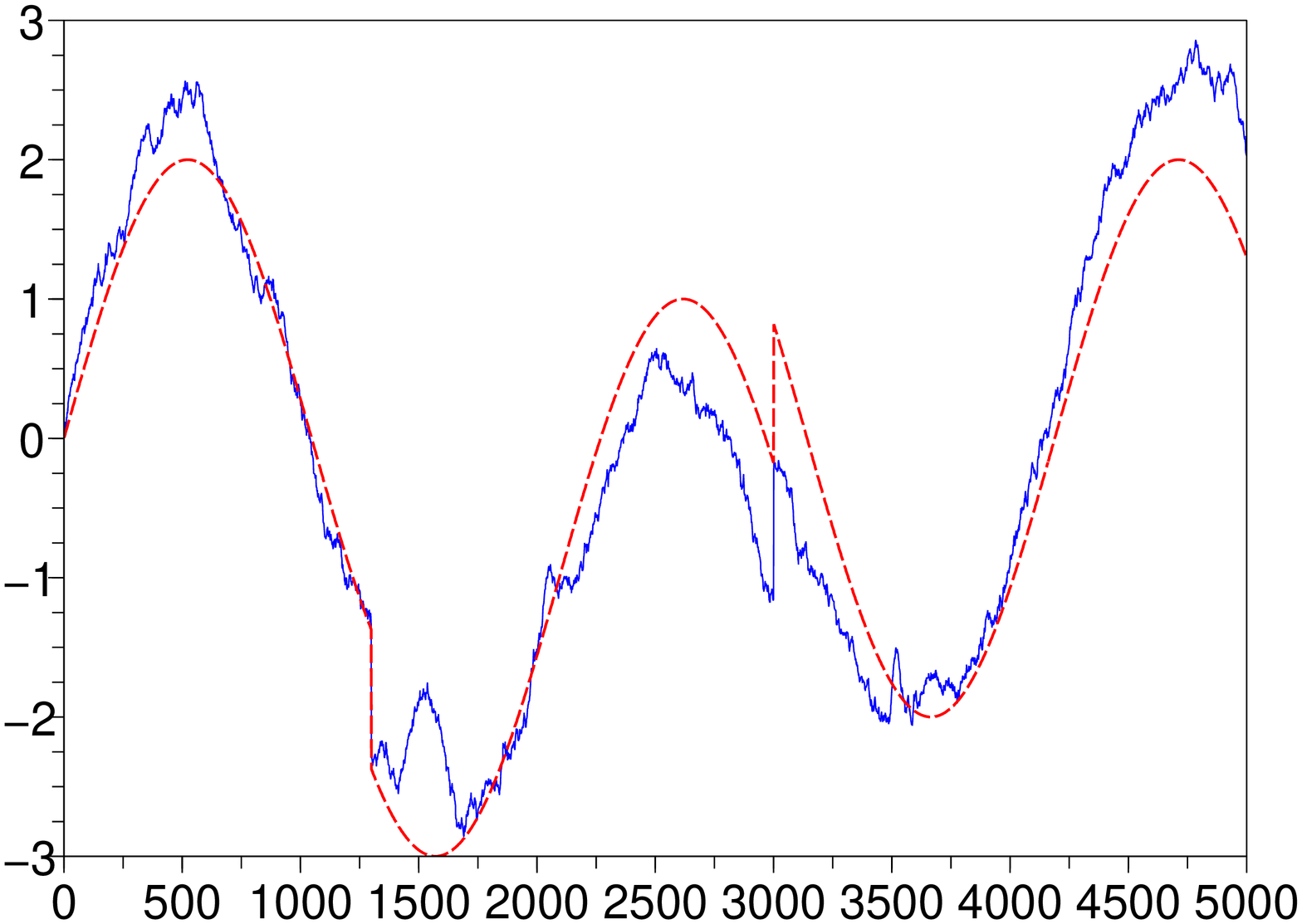}}}}
 {\subfigure[\footnotesize Change-point detection -- Exact (+)]{\rotatebox{-90}{\includegraphics[width=4cm]{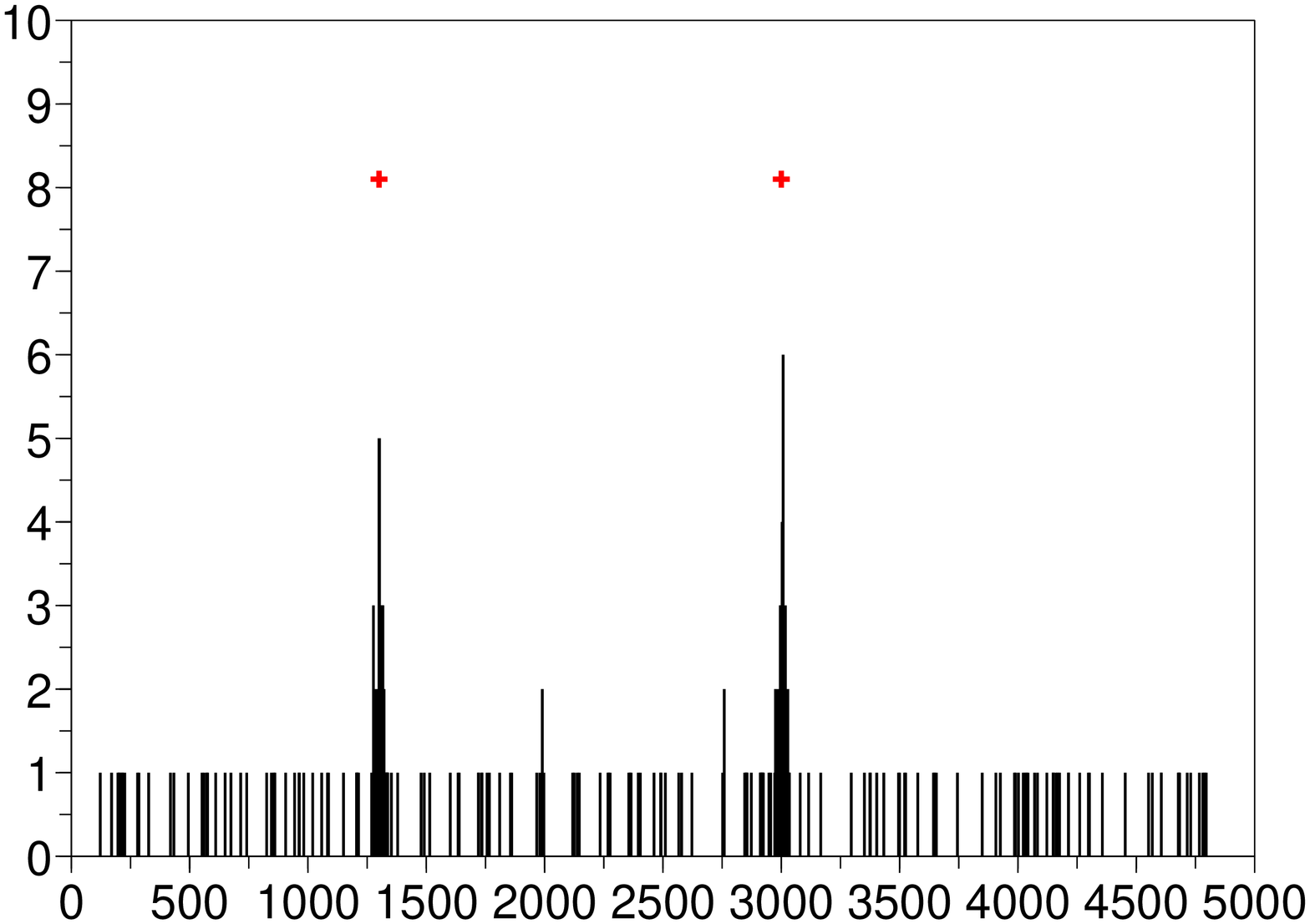}}}}
\caption{Sinusoidal signal -- Additive Perlin noise -- SNR: 10 db }
\label{Sp10}
\end{figure}

\begin{table*}\label{table}
\tiny\center

\begin{tabular}{|c|c|p{0.16cm}p{0.16cm}p{0.16cm}p{0.16cm}p{0.16cm}p{0.16cm}p{0.16cm}p{0.16cm}p{0.16cm}|c|}
\hline
 Noise type & SNR &\multicolumn{9}{c|}{Estimated number of segments} & Figure \\
&in&\multicolumn{9}{c|}{\bf True value in bold font}&references\\
&DB&1&2&3&4&5&6&7&8&$\geq$9&\\ \hhline{============}
Normal
(+)&0&0&0&0&1&{\bf 98}&1&0&0&0& figure \ref{M1}\\ \hline
Normal
(+)&-6&0&0&1&8&{\bf 79}&12&0&0&0&figure \ref{M2}\\  \hhline{============}

Normal
(+)&25&0&0&0&0&9&{\bf 83}&8&0&0& figure \ref{Bg25}\\ \hline Uniform
(+)&25&0&0&0&0&12&{\bf 81}&7&0&0& figure \ref{Bu25}\\ \hline Perlin
(+)&20&0&0&0&2&25&{\bf 40}&16&6&11& figure \ref{Bp20}\\ \hline
Uniform ($\times$)&20&0&0&0&1&11&{\bf 74}&14&0&0& figure
\ref{Bmu20}\\ \hhline{============} Normal
(+)&25&0&0&{\bf100}&0&0&0&0&0&0& figure \ref{Sg25}\\ \hline Normal
(+)&20&0&4&{\bf96}&0&0&0&0&0&0& figure \ref{Sg20}\\ \hline Perlin
(+)&10&0&18&{\bf42}&16&8&3&3&2&8& figure \ref{Sp10}\\ \hline
\end{tabular}

\caption{Summary of simulation results $\bullet$ +: additive noise
$\bullet$ $\times$: multiplicative noise \label{simuls}}
\end{table*}

%
%



\newpage

\begin{thebibliography}{3}
%
%

\bibitem{basseville}Basseville, M., Nikiforov, I.V.: Detection of
Abrupt Changes: Theory and Application. Prentice-Hall (1993).
Available online at {\tt http://www.irisa.fr/sisthem/kniga/}.

\bibitem{belkoura}Belkoura L.: Change point detection with application
to the identification of a switching process. In: El Jai A., Afifi
L., Zerrik E. (eds) Systems Theory: Modelling, Analysis and Control,
Internat. Conf. Fes (Morocco), pp. 409-415, Presses Universitaires
de Perpignan (2009). Available online at {\tt
http://hal.inria.fr/inria-00363679/en/}.

\bibitem{belkoura2}
Belkoura L., Richard J.-P., Fliess M.: Parameters estimation of
systems with delayed and structered entries. Automatica {\bf 45},
pp. 1117-1125 (2009).

\bibitem{brodsky1}Brodsky, B.E., Darkhovsky, B.S.: Nonparametric Methods in
Change-Point Problems. Kluwer (1993).

\bibitem{brodsky2}Brodsky, B.E., Darkhovsky, B.S.: Non-Parametric Statistical
Diagnosis: Problems and Methods. Kluwer (2000).

\bibitem{chambert} Chambert-Loir, A.: Alg\`{e}bre corporelle, \'Editions \'Ecole
Polytechnique (2005). English translation:  A Field Guide to
Algebra. Springer (2005).

\bibitem{chen}Chen, W.K.: Passive and Active Filters: Theory and
Implementations. Wiley (1986).

\bibitem{csorgo}Cs\"{o}rg\"{o}, M., Horv\'{a}th, L.: Limit Theorems in Change-Point
Analysis. Wiley (1997).

\bibitem{drag}Dragotti, P.L., Vitterli, M.: Wavelets footprints:
theory, algorithms, and applications. IEEE Trans. Signal Proc. {\bf
51}, pp. 1306-1323 (2003).

\bibitem{stanley}
Flajolet P., Sedgewick R.: Analytic Combinatorics. Cambridge
University Press (2009).

\bibitem{ans}Fliess, M.: Analyse non standard du bruit.
C.R. Acad. Sci. Paris Ser. I {\bf 342}, pp. 797-802 (2006).

\bibitem{arima}Fliess, M.: Critique du rapport
signal \`{a} bruit en communications num\'{e}riques. ARIMA \textbf{9}, pp.
419-429 (2008). Available online at {\small \tt
http://hal.inria.fr/inria-00311719/en/}.


\bibitem{finance}Fliess M., Join C.: Towards new technical indicators for trading systems and
risk management. \emph{15$^{th}$ IFAC Symp. System Identif.},
Saint-Malo (2009). Available online at {\tt
http://hal.inria.fr/inria-00370168/en/}.

\bibitem{finance2}Fliess M., Join C.: Systematic risk analysis:
first steps towards a new definition of beta. \emph{COGIS'09}, Paris
(2009). Available online at {\tt
http://hal.inria.fr/inria-00425077/en/}.

\bibitem{easy}Fliess, M., Join, C., Sira-Ram\'{\i}rez, H.: Non-linear estimation
is easy. Int. J. Modelling Identification Control. {\bf 4}, pp.
12-27 (2008).


\bibitem{gretsi}
Fliess, M., Join, C., Mboup, M., Sira-Ram\'{\i}rez, H.: Analyse et
repr\'esentation de signaux transitoires : application \`a  la
compression, au d\'ebruitage et \`a  la d\'etection de ruptures.
{\em $20^e$ Coll. GRETSI}, Louvain-la-Neuve (2005). Available online
at {\tt http://hal.inria.fr/inria-00001115/en/}.

\bibitem{mexico}Fliess, M., Mboup, M., Mounier, H., Sira-Ram\'{\i}rez, H.: Questioning some
paradigms of signal processing via concrete examples. In:
Sira-Ram\'{\i}rez H., Silva-Navarro G. (eds.) Algebraic Methods in
Flatness, Signal Processing and State Estimation, pp. 1-21,
Editiorial Lagares (2003). Available online at {\tt
http://hal.inria.fr/inria-00001059/en/}.

\bibitem{esaim}Fliess, M., Sira-Ram\'{\i}rez, H.: An algebraic
framework for linear identification, ESAIM Control Optim. Calc.
Variat. {\bf 9}, pp. 151-168 (2003).

\bibitem{garnier}Fliess, M., Sira-Ram\'{\i}rez, H.: Closed-loop parametric
identification for continuous-time linear systems via new algebraic
techniques. In: Garnier, H., Wang, L. (eds) Identification of
Continuous-Time Model Identification from Sampled Data, pp. 363-391,
Springer (2008).

\bibitem{gij}Gijbels, I., Hall, P., Kneip, A.: On the estimation of
jump points in smooth curves. Ann. Instit. Statistical Math. {\bf
51}, pp. 231-251 (1999).

\bibitem{lavielle}Lavielle, M.: Using penalized contrasts for change-point
problem. Signal Processing {\bf 85}, pp. 1501-1510 (2005).

\bibitem{lebarbier}Lebarbier, E.: Detecting mutiple change-points in the mean of
a Gaussian process by model selection. Signal Processing {\bf 85},
pp. 717-736 (2005).

\bibitem{lobry}Lobry, C., Sari, T.: Nonstandard analysis and representation
of reality. Int. J. Control \textbf{81}, pp. 517-534 (2008).

\bibitem{mallat}Mallat, S.: A Wavelet Tour of Signal Processing
(2$^{nd}$ ed.). Academic Press (1999).

\bibitem{mboup}Mboup M.: Parameter estimation for signals described by differential
equations. Applicable Anal.  \textbf{88}, pp. 29-52  (2009).

\bibitem{spike}Mboup M.: A Volterra filter for neuronal spike
detection. Preprint (2008). Available online at {\tt
http://hal.inria.fr/inria-00347048/en/}.

\bibitem{ajaccio}Mboup M., Join C., Fliess M.: A delay estimation
approach to change-point detection. {\em 16$^{th}$ Medit. Conf.
Control Automat.}, Ajaccio (2008). Available online at {\tt
http://hal.inria.fr/inria-00179775/en/}.

\bibitem{McC}McConnell, J., Robson, J.:  Noncommutative Noetherian Rings.
Amer. Math. Soc. (2000).

\bibitem{miku1}Mikusinski, J.:
Operational Calculus ($\text{2}^{nd}$ ed.), Vol. 1.  PWN \& Pergamon
(1983).

\bibitem{miku2}
Mikusinski, J., Boehme, T.: Operational Calculus ($\text{2}^{nd}$
ed.), Vol. 2. PWN \& Pergamon (1987).

\bibitem{ollivier}
Ollivier F., Moutaouakil S., Sadik B.: Une m\'{e}thode d'identification
pour un syst\`{e}me lin\'{e}aire \`{a} retards. C.R. Acad. Sci. Paris Ser. I
{\bf 344}, pp.709-714 (2007).

\bibitem{perlin}Perlin, K.: An image synthetizer.
ACM SIGGRAPH Comput. Graphics {\bf 19}, pp. 287-296 (1985).

\bibitem{pol}van der Pol, B., Bremmer, H.: Operational Calculus Based
on the Two-Sided Laplace Integral ($\text{2}^{nd}$ ed.). Cambridge
University Press (1955).

\bibitem{singer}van der Put, M., Singer, M.F.: Galois Theory of Linear
Differential Equations, Springer (2003).

\bibitem{sinica}Raimondo, M., Tajvidi, N.: A peaks over threshold model
for change-points detection by wavelets. Statistica Sinica {\bf 14},
pp. 395-412 (2004).

\bibitem{rudolph}
Rudolph J., Woittennek F.: Ein algebraischer Zugang zur
Parameteridentifkation in linearen unendlichdimensionalen Systemen.
at--Automatisierungstechnik {\bf 55}, pp. 457-467 (2007).

\bibitem{zoran}Tiganj Z., Mboup M.: Spike detection and sorting: combining
algebraic differentiations with ICA. {\em 8$^{th}$ Int. Conf. Indep.
Component Anal. Signal Separat.}, Paraty, Brazil (2009). Available
online at {\tt http://hal.inria.fr/inria-00430438/en/}.

\bibitem{tourneret}Tourneret, J.-Y., Doisy, M., Lavielle, M.:
Bayesian off-line detection of multiple change-points corrupted by
multiplicative noise: application to SAR image edge detection.
Signal Processing {\bf 83}, pp. 1871-1887 (2003).

\bibitem{yosida}Yosida, K.: Operational Calculus: A Theory of
Hyperfunctions ({\it translated from the Japanese}). Springer
(1984).

\end{thebibliography}


\end{document}